\documentclass[12pt]{article}
\usepackage{amsmath}
\usepackage{amssymb}
\usepackage{url}

\usepackage[T1]{fontenc}
\usepackage{mathptmx}
\usepackage{diagrams}
\diagramstyle{midshaft,PostScript=dvips,nohug}
\newarrow{Dotsto}.. ..>

\newtheorem{lemma}{Lemma}[section]
\newtheorem{prop}[lemma]{Proposition}
\newtheorem{thm}[lemma]{Theorem}

\newtheorem{defn}[lemma]{Definition}

\newcommand{\ev}{\operatorname{ev}}
\newcommand{\id}{\operatorname{id}}

\newcommand{\FM}{\mathfrak{M}}

\newcommand{\X}{\mathcal{X}}
\newcommand{\W}{\mathcal{V}}
\newcommand{\Z}{\mathcal{Z}}

\newcommand{\Y}{\mathcal{Y}}

\newcommand{\B }{\mathcal{B}}

\newcommand{\pit}{\pitchfork}

\newcommand{\of}{\overline{F}}
\newcommand{\uf}{\underline{F}}

\title{Bundle functors and fibrations}

\author{Anders Kock}
\date{}

\begin{document}

\maketitle

\section*{Introduction}The notions of bundle, and bundle functor, are   
useful and well exploited notions in topology and differential geometry, cf.\ e.g.\ 
\cite{[KMS]}, as well as in other branches of mathematics. The category 
theoretic set up relevant for these notions is that of fibred 
category, likewise a well exploited notion, but for certain considerations 
in the context of bundle functors, it can  be carried further. In particular, we formalize and develop, 
in terms of fibred categories, some 
of the differential geometric constructions: tangent- and cotangent bundles, 
(being examples of bundle functors, respectively star-bundle 
functors, as in \cite{[KMS]}), as well as 
jet bundles (where the formulation of the functorality properties,
 in terms of fibered categories, is probably 
new).

Part of the development in the present note was expounded in \cite{[DFET]}, 
and is repeated  almost verbatim in the Sections \ref{FSX} and 
\ref{DFX} below. These 
sections may have interest as a piece of pure category theory, not 
referring to differential geometry.

\section{Basics on Cartesian arrows}\label{BCAX} 
We recall here some classical notions. 

Let $\pi : \X \to \B$ be any functor. For $\alpha :A \to B$ in $\B$, 
and for objects $X, Y\in \X$ with $\pi (X)=A$ and $\pi (Y)=B$,
let $\hom_{\alpha}(X,Y)$ be the set of arrows $h:X\to Y$ in $\X$ with $\pi 
(h)=\alpha$.

The {\em 
fibre} over $A\in \B$ is the category, denoted $\X_{A}$, whose objects are the $X \in 
\X$ with $\pi (X)=A$, and whose arrows are arrows in $\X$ 
which by $\pi$ map to  $1_{A}$; such arrows are 
called {\em vertical} (over $A$). The hom functor of $\X _{A}$ is 
denoted $\hom _{A}$.

Let $h$ be an arrow $X \to Y$, and denote  $\pi (h)$ by  $\alpha 
: A \to B$, where $A=\pi(X)$ and $B=\pi (Y)$.  
 For any arrow $\xi : C\to A$, and any object $Z\in \X_{C}$, 
post-composition with $h:X \to Y$ defines a map $$h_{*}:\hom 
_{\xi}(Z,X)\to \hom_{\xi . \alpha} (Z,Y).$$
(we compose from left to right). Recall that $h$ is called {\em 
Cartesian} (with respect to $\pi$) if this map is a bijection, for all such $\xi$ and $Z$.
It is easy to see that Cartesian arrows compose, and that 
isomorphisms are Cartesian. In particular, the Cartesian arrows form 
a subcategory of $\X$. Also, a vertical Cartesian arrow is an 
isomorphism.

\medskip

\noindent{\bf Example.} The following is a fundamental example, which 
will also be the origin for some of the applications that we present. Let $\B$ be 
any category, and let $\B^{2}$ be the category of arrows in $\B$, so 
the arrows in $\B^{2}$ are the commutative squares in $\B$. Let 
$\partial_{1}:\B^{2}\to \B$ be the functor which to an arrow assigns 
its codomain. Then a commutative square in $\B$ is a Cartesian 
arrow in $\B^{2}$ 
(with respect to $\partial_{1}: \B^{2}\to \B$) precisely when the 
square is a pull-back.

\medskip

The property of being  a Cartesian arrow is clearly a kind of universal 
property.
There is a weaker notion of when an arrow $h$ as above is {\em pre-Cartesian}
\footnote{called {\em Cartesian} in the earlier literature 
(Grothendieck et.al); we follow the 
terminology mostly in use now, see e.g.\ \cite{[Borceux]}}, namely that, 
for any $Z\in \X_{A}$, post-composition with $h$
 defines a bijection $$h_{*}:\hom 
_{A}(Z,X)\to \hom_{ \alpha} (Z,Y).$$
The property of being  a pre-Cartesian arrow is clearly a  universal 
property. In  fact, to say that $h:X \to Y$ is preCartesian over $\alpha$ can be 
expressed by saying that $h$ is terminal in a certain ``relative 
comma-category'' $\X_{A}\downarrow_{\alpha}Y$ whose objects are 
arrows in $\X$  over $\alpha$ with codomain $Y$, and whose arrows are 
arrows in $\X_{A}$ making an obvious triangle commute.

\medskip

\noindent{\bf Remark.} There are dual notions of coCartesian and pre-coCartesian arrows; 
they will not play much role in the following, except that we at one 
point shall consider the latter notion; thus, if $\alpha: A \to B$ in $\X$, 
and $X\in \X_{A}$, we have another ``relative comma-category'' $X 
\downarrow _{ \alpha} \X_{B}$ whose objects are 
arrows in $\X$  over $\alpha$ with domain $X$, and whose arrows are 
arrows in $\X_{B}$ making an obvious triangle commute. Then a 
pre-coCartesian arrow over $\alpha$ with domain $X$ is by definition 
an initial object in this category.

\medskip

Clearly, if $h$ is Cartesian, then it is also 
pre-Cartesian. Also Cartesian arrows over $\alpha$, with given 
codomain $Y$, are unique up to unique vertical map (necessarily 
invertible) in 
$\X_{A}$; the same applies to pre-Cartesian arrows.

If $h$ is Cartesian, the injectivity of $h_{*}$ implies the 
cancellation property that $h$ is ``monic w.r.\ to $\pi$'', 
meaning that  for parallel arrows $k,k'$ in $\X$ with 
codomain $X$, and with $\pi (k)=\pi (k')$, we have that 
$k.h = k'.h \mbox { implies } k=k'$. 

For later use, we recall a basic fact: 
\begin{lemma}\label{lemmax} If $k=k'.h$ is Cartesian and 
$h$ is Cartesian, then $k'$ is Cartesian.\end{lemma}

 For pull-back squares, this is well known property  for the functor 
$\partial _{1}: \B^{2}\to \B$, cf.\ the Example above.
  
\medskip

The functor $\pi :\X \to \B$ is called a {\em fibration} if there are 
enough Cartesian arrows, in the following sense: for every 
$\alpha : A \to B$ in $\B$ and every $Y \in \X_{B}$, there 
exists a Cartesian arrow over $\alpha$ with codomain $Y$. Such arrow 
is called {\em a Cartesian lift} of $\alpha$ with codomain $Y$.  
A choice, for each arrow $\alpha :A \to B$ in $\B$ 
and for each $Y\in \X_{B}$, of a Cartesian lift of $\alpha$  with 
codomain $Y$, is called a {\em clevage} of the fibration 
$\pi$. The domain 
of the chosen Cartesian arrow over $\alpha$ with codomain $Y$ is sometimes 
denoted $\alpha^{*}(Y)$. We use clevages mainly as a notational tool, to facilitate 
reading, but generally, we avoid clevages.

The functor $\partial_{1}:\B^{2}\to \B$ (cf.\ the Example above) is a fibration precisely when 
$\B$ is a category with pull-backs; then $\partial _{1}$ is called ``the codomain fibration''. 
A cleavage for it amounts to a choice of pull-back diagrams in $\B$.

\section {The ``factorization system'' for a fibration}\label{FSX}
Let $\pi: \X \to \B$ be a fibration, and let $z:Z\to Y$ be an arrow in $\X$. 
Let $h:X \to Y$ be a Cartesian arrow over $\alpha:= \pi (z)$. By the 
universal property of the Cartesian arrow $h$, there is a unique 
vertical $v:Z \to X$ with $z=v.h$. 

Thus, every arrow $z$ in $\X$ may be written as a composite 
of a vertical arrow followed by a Cartesian arrow (Cartesian arrows, we like to 
think of as being ``horizontal''). And, crucially, 
this decomposition of $z$ is unique modulo a {\em unique} vertical
isomorphism. Or, equivalently, modulo a unique arrow which is at the same 
time vertical and cartesian. (Recall that for vertical arrows, 
being Cartesian is equivalent to being an isomorphism (= invertible).) This means that every arrow $z$ in $\X$ 
may be represented by a pair $(v,h)$ of arrows with $v$ vertical and 
$h$ cartesian, with $z=v.h$, (and with  the codomain of $v$ equal to the domain 
of $h$). We call such a pair a ``vh composition pair'', to make the 
analogy with vh spans, to be considered below, more explicit. Two such pairs $(v,h)$ and 
$(v',h')$ represent the same arrow in $\X$ iff there exists a vertical 
cartesian (necessarily unique, and necessarily invertible) $s$ such that 
\begin{equation}\label{eqx1}v.s=v' \mbox{ and } s.h'=h. \end{equation}
We say that $(v,h)$ and $(v',h')$ are {\em equivalent} if this holds.
The composition of arrows in $\X$ can be described in terms of 
representative vh composition pairs, as follows. If $z_{j}$ is represented by 
$(v_{j},h_{j})$ for $j=1, 2$, then $z_{1}.z_{2}$ is represented by 
$(v_{1}.w, k.h_{2})$, where $k$ is cartesian over $\pi (h_{1})$ and $w$ 
is vertical, and the square displayed commutes:
\begin{equation}\label{interpolx}\begin{diagram}\cdot&&&&&\\
\dTo^{v_{1}}&&&&\\
\cdot & \rTo ^{h_{1}}&\cdot&&&\\
\dTo^{w}&&\dTo_{v_{2}}&&\\
\cdot&\rTo_{k}&\cdot & \rTo_{h_{2}}&\cdot
\end{diagram}\end{equation}
Such $k$ and $w$ exists (uniquely, up to unique vertical cartesian 
arrows): 
construct first $k$ as a Cartesian lift of $\pi (h_{1})$ with same 
codomain as $v_{2}$, then use 
the universal property of Cartesian arrows to construct $w$.

The arrows $z_{1}$ and $z_{2}$ may be inserted, completing the 
diagram with two commutative triangles, since $z_{j}=v_{j}.h_{j}$. 
But if we refrain from doing so, we have a blueprint for 
 a succinct 
and choice-free description of the fibrewise dual $\X ^{*}$ of the fibration 
$\X\to \B$, to be described in Section \ref{DFX}. 

Note that a vh  factorization of an arrow in $\X$ is much reminiscent 
of the  factorization for an $E$-$M$ factorization system, as in 
\cite{[Borceux]} I.5.5, say, (with the 
class of  
vertical arrows playing the role of $E$, and the class of Cartesian 
arrows  playing the 
role of $M$; however, note that not every isomorphism in $\X$ is 
vertical.

\section{Construction of functors out of a fibered category}\label{CFX}
Let $\pi :\X \to \B$ be a fibration, 
and consider a functor $F: \X \to \Y$. Let $\underline{\X}$ denotes the 
category of vertical arrows in $\X$.
Then by restriction, $F$ gives 
a functor  $\underline{F}: \underline{\X} \to \Y$. The restriction of 
$F$  (or of $\uf$) to the fibre $\X_{A}$ is denoted $F_{A}$. 
Similarly, let $\overline{\X}$ denote the 
category over $\B$ consisting of the Cartesian arrows of $\X$ only, 
and let $\of$ denote the restriction of $F$ to $\overline{\X}$. Then $F$ gives rise to the 
following data: 

1) for each $A\in \B$, a functor $F_{A}: \X_{A}\to \Y$

2 a functor   $\overline{F}: \overline{\X} \to \Y$.
\newline Since the $F_{A}$s and $\overline{F}$ are restrictions of the same 
functor $F$, it is clear that we have the properties

3) if $s$ is vertical over $A$, and Cartesian, $F_{A}(s) = 
\overline{F}(s)$

4) Given a commutative square in $\X$, with $v$ and $w$ vertical and 
with $h$ and 
$k$ Cartesian (over $\alpha :A \to B$, say)
$$\begin{diagram}\cdot & \rTo^{k}& \cdot\\
\dTo^{v}&&\dTo_{w}\\
\cdot & \rTo_{h} & \cdot
\end{diagram} $$
Then $ F_{A} (v).\of (h)= \of (k).F_{B}(w)$.

\begin{thm}\label{constrx}Given functors $F_{A}:\X_{A}\to \Y$ (for all $A \in \B$), 
and given a functor $\of : \overline{\X} \to \Y$ over $\B$ as in 1) 
and 2), and assume that the conditions 3) and 4) hold. Then there 
exists a unique functor $F: \X \to \Y$  with restrictions 
$F_{A}$ to the fibres $\X_{A}$ and with restriction $\of$ to the Cartesian 
arrows. 

If $\Y$ comes with a functor to  $\B$, and if  for all $A\in \B$, 
$F_{A}$ factors through $\Y_{A} \subseteq \Y$, then the constructed 
$F$ is a functor over $\B$.
\end{thm}
{\bf Proof.} Given an arrow $x$ over $\alpha :A \to B$, say. Since 
$\X \to \B$ is a fibration, $x$ admits a vh factorization $x=v.h$ with 
$v$ vertical over $A$ and $h$ Cartesian over $\alpha$, so we are 
forced to define $F(x):= F_{A}(v).\of (h)$. To see that this $F$ is 
well defined, we consider another possible vh factorization $x=v'. 
h'$. It compares with the given $v.h$ by a vertical Cartesian $s$ 
with $v'= v.s$ and $h=s.h'$. We have
$$F_{A}(v). \of (h)= F_{A}(v).\of (s.h') = F_{A}(v). \of (s).\of 
(h'),$$
$$F_{A}(v').\of (h')=F_{A}(v.s).\of (h') =F_{A}(v).F_{A}(s).\of 
(h'),$$
using that $\of$ and $F_{A}$ are functors. By condition 3), $\of (s)= 
F_{A}(s)$, so the two expressions agree.
Let us next prove that the $F$ constructed preserves composition 
of arrows, say $f_{1}.f_{2}$. Pick a vh factorization of $f_{1}$, 
say $f_{1}=v_{1}.h_{1}$, 
and similarly for $f_{2}$. Interpolate a $w$ and $k$ as in 
(\ref{interpolx}); then use condition 4) for the interpolated square.
--The last assertion is obvious.

\medskip

Even when, as in the last statement of the Proposition, the category 
$\Y$ is given as a category over $\B$, it is not assumed to be 
fibered over $\B$. But if $\Y \to \B$ happens to be a fibration, then 
given the family of functors $F_{A}:\X_{A}\to \Y_{A}$, the data of 
the functor $\of$  may be formulated in an alternative way, provided 
we assume given cleavages of both $\X\to \B$ and $\Y \to \B$. For 
then, to give the value of $\of$ on a Cartesian arrow $X'\to X$ over 
$\alpha$, it suffices to give the value of $\of$ on the chosen 
Cartesian arrow  $h: \alpha^{*}(X) \to X$ over $\alpha$ and with codomain 
$X$. 
This value is an arrow in $\Y$ over $\alpha$ and with codomain $F_{B}(X)$, 
and as such factors uniquely by a vertical arrow $v_{\alpha, X}$
followed by the chosen Cartesian arrow over $\alpha$ with that 
codomain, thus
$$\begin{diagram}F_{A}(\alpha^{*}(X))&&&\\
\dTo^{v_{\alpha, X}}&\rdTo ^{\of (h)}&\\
\alpha^{*}(F_{B}(X))&\rTo&F_{B}(X)
\end{diagram}$$
with the bottom arrow the chosen Cartesian. (Note that $\of (h)$ 
need not be Cartesian; we did not assume that $\of$ preserves the 
property of being Cartesian.) So the data of $\of$ resides in the 
$F_{A}$s, together with the vertical maps
\begin{equation}\label{fibstrx}v_{\alpha,X}:F_{A}(\alpha^{*}(X))\to 
\alpha^{*}(F_{B}(X)). 
\end{equation}
(The $v_{\alpha, X}$ thus derived satisfy certain equations, in 
particular, for fixed $\alpha$, $v_{\alpha ,X}$ is natural in $X \in 
\X_{B}$; there are also equational conditions involving the comparison 
isomorphisms between $\alpha^{*}\circ \beta^{*}$ and $(\beta \circ 
\alpha )^{*}$.
In terms of pseudofunctors sometimes used to present fibrations,  
$v$ is a lax (or colax?) transformation between the pseudofunctors 
representing $\X$ and $\Y$, respectively. - We shall not enter into 
these conditions, since the conditions in the Theorem are clear 
enough.)

\section{The dual fibration $\X^{*}$; comorphisms}\label{DFX}

The construction\footnote{The 
construction,  in elementary terms, 
of the dual fibration
 can be distilled 
out of a more general construction \cite{[BGN]} by Barwick et al.\ in the context 
of quasi-categories. I was unaware of their construction when I put a 
preliminary version \cite{[DFET]} of the present paper on arXiv. I want to thank 
them for calling my attention to their work. The  construction (for 
categories, not for quasi-categories) was 
apparently also known by Borceux, cf.\ Exercise 8.8.2 in \cite{[Borceux]} II.} 
presented in this Section is  elementary. (In fact it is clear that it 
makes sense for categories and fibrations internal to an exact 
category.) 
It is is a direct  
generalization of the ``star bundle'' construction of \cite{[KMS]} 41.1, 
where it is presented to 
account for the functorial properties of, say,  the formation of cotangent 
bundles in differential 
geometry.

\medskip

Given a fibration $\pi: \X \to \B$. We describe another category 
$\X^{*}$ over $\B$, the ``fibrewise dual of $\X \to \B$'', as follows: The objects of $\X^{*}$ are the same as those of $\X$; 
the arrows $X \to Y$ are represented by vh spans, in the following 
sense: 
\begin{defn}A {\em vh span} in $\X$ from $X$ to $Y$ is a diagram in $\X$ of the form
\begin{equation}\label{vhx}\begin{diagram}\cdot & \rTo^{h}& Y\\
\dTo^{v}&&\\
X&&
\end{diagram}\end{equation}
with $v$ vertical and $h$  cartesian.
\end{defn}
The set of arrows in $\X^{*}$ from $X$ 
to $Y$ are equivalence classes of vh spans from $X$ to $Y$, for the 
equivalence relation $\equiv$ given by
$(v,h) \equiv (v',h')$ if there exists a vertical isomorphism $s$ 
(necessarily unique) in 
$\X$ so that 
\begin{equation}\label{eqx2}s.v.=v' \mbox{ and } s.h=h'. \end{equation}
We denote the equivalence class of the vh span $(v,h)$ by $\{(v,h)\}$. 
They are the arrows of $\X^{*}$;   the direction of a the arrow $\{(v,h)\}$ 
is determined by its 
Cartesian part $h$.
  
 Composition  has to be described in terms of  representative 
pairs; it is in fact the standard composite of spans, but let us be 
explicit: If $z_{j}$ is represented by 
$(v_{j},h_{j})$ for $j=1, 2$, then $z_{1}.z_{2}$ is represented by 
$(w.v_{1},k.h_{2})$, where $k$ is Cartesian over $\pi (h_{1})$ and $w$ is 
vertical, and the square displayed commutes:
\begin{equation}\label{compx}\begin{diagram}\cdot &\rTo^{k}&\cdot & \rTo^{h_{2}}&\cdot\\
\dTo^{w}&&\dTo_{v_{2}}&&\\
\cdot&\rTo_{h_{1}}&\cdot &&\\
\dTo^{v_{1}}&&&&\\
\cdot
\end{diagram}\end{equation}
Such $k$ and $w$ exists (uniquely, up to unique vertical cartesian 
arrows): 
construct first $k$ as a cartesian lift of $\pi (h_{1})$, then use 
the universal property of cartesian arrows to construct $w$. (The 
square displayed will then actually be a pull-back diagram, thus the 
composition described will be the standard composition of spans.)

Composition of vh spans  does not give a definite vh span, 
but rather an 
equivalence class of vh spans. So referring to (\ref{compx}), the 
composite of  $\{(v_{1},h_{1})\}$ with  $\{(v_{2},h_{2})\}$ is defined by
$$   \{(v_{1},h_{1})\}.\{(v_{2},h_{2})\}:= \{(w.v_{1},k.h_{2})\}.$$

There is a functor $\pi^{*}$ from $\X^{*}$ to $\B$; on objects, it agrees with $\pi 
: \X \to \B$; and   $\pi ^{*}(\{(v,h)\})= 
\pi (h)$. Note that if 
$v:X'\to X$ is vertical, the vh span $(v,1)$ represents a morphism 
$X \to X'$ in $\X^{*}$.  

Clearly, a vertical arrow in $\X^{*}$ has a unique representative 
span of the form $(v,1)$.
So  the fibres of $\pi^{*}:\X^{*} \to \B$ are 
canonically isomorphic
to the 
duals of the fibres of $\pi: \X \to \B$, i.e.\ $(\X ^{*})_{A} \cong (\X 
_{A})^{op}$; so 
 $\X^{*}$ is ``fibrewise dual'' to $\X$ (but is not 
in general dual to $\X$, since the functor $\pi^{*}:\X^{*}\to \B$ is 
still a covariant functor). 
The arrows in $\X^{*}$, we call {\em comorphisms} in $\X$; it is ususally 
harmless to use the name ``comorphism'' also for a representing vh span $(v,h)$.

There are two special classes of comorphisms: the first class 
consists of those comorphisms that 
can be represented by a pair $(v,1)$ where $1$ is the relevant identity arrow. They 
are precisely the vertical arrows for $\X^{*}\to \B$.  -- The second class 
consists of those comorphisms that 
can be represented
  by a pair $(1,h)$ where $1$ is 
the relevant identity arrow. We shall see that these are precisely 
the cartesian morphisms in $\X^{*}$. 

We first note that if $(v,h)$ represents an arbitrary arrow in $\X^{*}$, 
then 
\begin{equation}\label{factx}(v,h) \in \{(v,1)\}.\{(1,h)\} ;
\end{equation}
this is witnessed by the diagram
$$
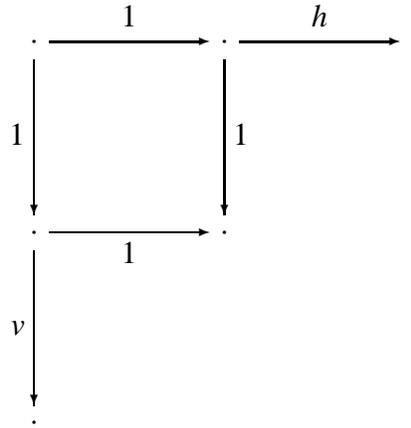
\begin{diagram}\cdot& \rTo ^{1}&\cdot & \rTo^{h}&\cdot\\
\dTo^{1}&&\dTo_{1}&&\\
\cdot&\rTo_{1}&\cdot&&\\
\dTo^{v}&&&&\\
\cdot&&&&
\end{diagram}$$
since the upper left square is of the form considered in 
(\ref{compx}).

\begin{prop} An arrow $g$ is Cartesian in $\X^{*}$ iff it admits a vh 
representative of the form $(1,h)$. Any vh representative of such $g$ 
is of the form $(w,k)$ with $w$ (vertical and) invertible.
\end{prop}
{\bf Proof.} In one direction, let $(1,h)$ represent a comorphism 
$Y\to Z$ over the arrow $\beta$ in $\B$, and let 
$(v, k)$ represent  a comorphism $X\to Z$  over $\alpha. \beta$. We 
display these data as the full arrows in the following display (in 
$\X$ and $\B$):
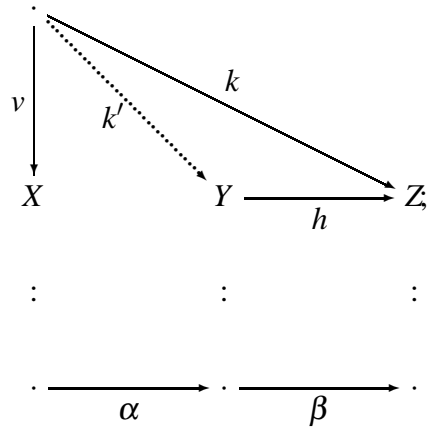
$$ \begin{diagram}\cdot&&&&\\
\dTo^{v}&\rdDotsto_{k'} \rdTo(4,2)^{k}&&&\\
X&&Y & \rTo_{h}& Z\\
:&&: && :\\
\cdot&\rTo_{\alpha }&\cdot &\rTo _{\beta }&\cdot 
\end{diagram};$$
The dotted arrow  $k'$, with $k'.h=k$, comes about by using the universal property of the 
Cartesian arrow $h$ in $\X$. Since $k$ and $h$ are Cartesian, then so 
is $k'$, by  Lemma \ref{lemmax}.
 So $(v,k')$ is a comorphism 
over $\alpha$, and $(v,k').(1.h) \equiv (v,k)$, and using the cancellation 
property of Cartesian arrows,  $(v,k')$ is easily 
seen to represent the unique comorphism over $\alpha . \beta$ composing with 
$(1,h)$ to give  $(v,k)$.  

In the other direction, let $g$ be a cartesian arrow in $\X^{*}$. 
Let $(w,k)$ be an arbitrary representative of $g$. Then 
by (\ref{factx}), $g=\{(w,1)\}.\{(1,k)\}$. Since $g$ is assumed cartesian in 
$\X^{*}$, and $\{(1,k)\}$ is cartesian  by what is already proved, it 
follows from Lemma \ref{lemmax} 
that $\{(w,1)\}$ is 
cartesian. Since it is also vertical, it follows that it is an 
isomorphism in $\X^{*}$, hence $w$ is an isomorphism in $\X$. (And 
this proves the second assertion of the Proposition.) Since 
$k$ is cartesian in $\X$, $w^{-1}.k$ is cartesian as well, and
$$(w,k)\equiv (1,w^{-1}.k),$$
so $g$ has a representative of the claimed form.

\begin{prop}The functor $\pi^{*}:\X^{*}\to \B$ is a fibration over $\B$
\end{prop}
{\bf Proof.} Let $\beta :A \to B$ be an arrow in $\B$, and let $Y\in 
\X_{B}$. Since $\X \to \B$ is a fibration, there exists in 
$\X$ a cartesian arrow $h$ over $\beta$ with codomain $Y$, 
and then the vh span $(1,h)$ 
represents, by the above, a cartesian arrow in $\X^{*}$ over $\beta$.

\medskip
The argument gives what may  briefly be expressed: the Cartesian arrows of $\X$ are the same as the 
Cartesian arrows of $\X^{*}$.

\medskip

Since $\X^{*}\to \B$ is a fibration, we may ask for its fibrewise 
dual $\X ^{**}$: 
\begin{prop}There is a canonical isomorphism over $\B$ between $\X$ 
and $\X^{**}$.
\end{prop}
{\bf Proof.} We describe an explicit functor $y: \X \to \X^{**}$.
Let us denote arrows in $\X^{*}$ by dotted arrows; they may be 
presented by vh spans $(v,h)$ in $\X$.
We first describe $y$ on vertical and cartesian arrows separately. For a 
vertical $v$ in $\X$, say $v:X\to X'$, we have the  
vh span $(v,1)$ in $\X$, which represents     a 
vertical 
arrow $\overline{v}: X'\dasharrow X$ in $\X^{*}$; thus we have a 
vh span $(\overline{v},1)$ in $\X^{*}$, which in turn represents
 a vertical arrow $X\to X'$ in $\X^{**}$. This arrow, we take as $y(v) \in \X^{**}$. 
Briefly, $y(v)= ((v,1),1)$. -- For a cartesian $h:X'\to Y$ (over 
$\beta$, say), we have a vh span $(1, h)$ in $\X$, which represents  a horizontal  arrow $\overline{h}: X'\dasharrow Y$ in $\X^{*}$ (cartesian 
over $\beta$); 
thus we have a vh span $(1,\overline{h})$ in $\X^{*}$, hence an 
arrow in $\X^{**}$, from $X'$ to $Y$ which we take as $y(h) \in \X^{**}$;
 briefly, $y(h)= (1,(1,v))$. The construction Theorem \ref{constrx} 
can now be applied; thus for 
a general $f:X\to Y$ in $\X$, we factor it $v.h$ with $v$ 
vertical and $h$ cartesian, and put $y(f):= y(v).y(h)$. We leave to 
the reader to verify the conditions  3) and 4) of the Theorem, i.e.\ that a different choice of $v$ and $h$ gives an equivalent 
vh span in $\X^{*}$, thus the same arrow in $\X^{**}$. 

Conversely, given an arrow $g: X \to Y$ in $\X^{**}$, represent it by 
a vh span  
in $\X^{*}$, $(\overline{v}, \overline{h})$,
$$\begin{diagram}X'& \rDotsto ^{\overline{h}}&Y\\
\dDotsto ^{\overline{v}}&&\\
X&&
\end{diagram}$$
Since $\overline{v}$ is vertical, we may pick a 
representative of $\overline{v}$ in the form $(v,1)$ with $v:X\to X'$,
 and since $\overline{h}$ is cartesian in $\X^{*}$, 
we may pick a representative of it if the form $(1,h)$, with $h:X'\to Y$ 
in $\X$. Then the composite $v.h : X\to Y$ makes sense in $\X$, and it goes by 
$y$ to the given $g$.

\medskip

For simplicity of notation and reading, one sometimes assumes that 
one has a 
cleavage for a given fibration $\X \to \B$, i.e.\  a choice of Cartesian arrows; for $\alpha$ an arrow 
in $\B$ and $Y$ an object in $\X$ over the codomain of $\alpha$,  
 the  chosen Cartesian arrow over $\alpha$ is denoted  $\alpha^{*}(Y) 
\to Y$. With such a cleavage, each equivalence class of vh 
composition pairs has a unique representative with one of these 
chosen arrows as h-part; and similarly for vh spans.

\section{The codomain fibration, and bundle functors}\label{CDFX}

Recall that if $\B$ is a category with pull-backs, then 
$\partial_{1}:\B^{2}\to \B$ is a fibration; the Cartesian arrows are 
the pull-back squares. This fibration is called the {\em codomain fibration}.
Note that for $A\in \B$, the 
category $(\B^{2})_{A}$ is the slice category $\B/A$.

 For simplicity 
of notation, we assume in this Section a cleavage, which here 
amounts to a choice of pull-back dia\-grams, for any $\alpha$ and $y$ 
with common co\-domain; then the following uses of the notation 
$\alpha^{*}$ is standard:
$$\begin{diagram}\cdot & \rTo^{}& \cdot\\
\dTo ^{ \alpha^{*}(y)}&&\dTo_{y}\\
A&\rTo_{\alpha}&B
\end{diagram} \mbox{\quad \quad or with slight abuse: \quad}\begin{diagram}\alpha^{*}(Y) & \rTo^{}& Y\\
\dTo ^{}&&\dTo_{y}\\
A&\rTo_{\alpha}&B.
\end{diagram} $$
Note that we do not assume that $\alpha^{*}(Y)= y^{*}(A)$.

The identity functor $\id_{B}: \B \to \B$ is likewise a 
fibration over $\B$ (this does not depend 
on $\B$ having pull-backs). When viewing $\B$ as being fibered over 
$\B$ in this way, it is sometimes useful to denote it $\B^{1}$, in 
analogy with $\B^{2}$; all arrows in $\B^{1}$ are Cartesian.

An important class of functors over $\B$ are  functors $T: \B^{1}\to 
\B^{2}$. Thus, the data of such functor amounts to a functor 
$T_{0}:\B\to \B$ plus a natural transformation $T : T_{0}\to 
\id_{\B}$. Such data is called a {\em bundle functor} in \cite{[KMS]} when 
$\B$ is the category of smooth manifolds (they require the instances 
of $\pi$ to be submersions).
Often, one does not notationally  distinguish between $T_{0}$ and 
$T$, or one writes $T$ for $T_{0}$ and $\pi$ for the natural 
transformation. 
 An example is the tangent bundle 
formation: If $A$ is a smooth manifold, $T(A)$ is the tangent bundle 
of $A$,  $\pi_{A}: T_{0}(A) \to A$ (we are disregarding for the moment the 
fibrewise vector space structure of the tangent bundle). Naturality of $\pi$ says that for 
$\alpha : A \to B$
$$\begin{diagram} T_{0}(A)& \rTo^{T_{0}(\alpha )}& T_{0}(B)\\
\dTo^{\pi_{A}}&&\dTo_{\pi_{B}}\\
A&\rTo_{\alpha}&B 
\end{diagram}$$
commutes. The bundle functor thus described is Cartesian (i.e.\ 
preserves Cartesian arrows) precisely 
when all squares of this form are pull-backs. (This square is clearly 
not a pull-back when  $T$ is the tangent bundle formation, unless 
$\alpha$ is a local diffeomorphism.)

\section{The fibrewise dual of the codomain fibration}\label{FDX}

We describe $(\B^{2})^{*}$, specializing the description in the 
previous Section. 
 Explicitly, 
for this special case, its objects 
are likewise arrows in $\B$, and the arrows over $\alpha :A \to B$, 
from $x:X\to A$ to 
$y: Y \to B$,  may be 
presented  in the form of  commutative diagrams $\eta$ 
(``comorphisms'' from $x$ to $y$)
\begin{equation}\label{comorx}\begin{diagram}\alpha^{*}(Y)&\rTo &Y\\
\dTo^{v}& &\\
X& \eta &\\
\dTo^{x} && \dTo_{y} \\
A&\rTo_{\alpha}&B
\end{diagram}\end{equation}
where the rectangle is a pull-back; since we have chosen pull-backs, the 
presentation is unique if we insist that the top arrow is a chosen 
Cartesian, as suggested by the notation. Note that, given the  $x$ and $y$, as well as and the 
$\alpha$, the information of the comorphism $\eta$ resides in 
the map denoted $v$. 

This kind of pull-back diagram was also considered in \cite{[Weber]}, under the name 
of ``pull-back around $\alpha$, $x$''. By the general theory of 
Section \ref{DFX}, the comorphism $\eta: x \dasharrow y$ exhibited in 
(\ref{comorx}) is Cartesian in $(\B ^2)^{*}$ iff $v$ is an 
isomorphism. This implies that it is a terminal object in the 
relative commacategory $(\B ^2)^{*}_{A} \downarrow _{\alpha}y$. We may 
also ask the dual question: when is $\eta$ pre-coCartesian, i.e.\ 
initial in the relative commacategory
$x\downarrow_{\alpha} ((\B ^2)^{*})_{B}$ ? This is precisely to say that 
the diagram is a {\em distributivity} pull-back, in the sense of 
\cite{[Weber]}; for, this means by definition  that it is a terminal object 
in the category of ``pull-backs around $\alpha$, $x$''  . The reason why our 
``initial'' then is substituted for  ``terminal'' in \cite{[Weber]} is just that, in 
our set up, the $A$-fibre of $(\B^{2})^{*}$ is {\em dual} to $\B/A$.  When 
the functor $\alpha^{*}: \B/B \to \B/A$ has a right adjoint 
$\Pi_{\alpha}$, then the $y$ occurring in (\ref{comorx}) is 
$\Pi_{\alpha}(x)$, and the $v$ is the back adjunction 
$\alpha^{*}\Pi_{\alpha }(x)\to x$.

\medskip

\noindent{\bf Remark.} It is worthwhile to reformulate the 
description of the fibrations $\B^{2}\to \B$ and $(\B^{2})^{*}\to \B$ 
for the case where $\B$ is the category of sets, so that an object 
$\xi :X\to A$ in $\B^{2}$ or in $(\B^{2})^{*}$ may be seen as a 
family $\{X_{a}\mid a\in A\}$ of sets $X_{a}:=\xi^{-1}(a)$. 

In this 
case, a morphism in $\B^{2}$ over $\alpha :A\to B$, from $X\to A$ to 
$Y\to B$, may be seen as a family of maps $\{ f_{a}:X_{a}\to 
Y_{\alpha (a)}\mid a\in A\}$, and a morphism in $(\B^{2})^{*}$ (i.e.\ a 
comorphism) over $\alpha:A\to 
B$ from $X\to A$ to $Y\to B$ may be seen as a family of maps
$\{f_{a}:Y_{\alpha(a)}\to X_{a}\mid a\in A\}$. Let us write $f: X 
\dasharrow 
Y$ for such a comorphism, reserving the plain arrow for acual set 
maps. Composition of comorphisms is 
essentially just composition of maps: if $f:X \dasharrow  Y$, as above, is a 
comorphism over $\alpha :A \to B$ and $g:Y\dasharrow Z$ similarly is a 
comorphism over $\beta: B\to C$, the composite of $f$ followed by  
$g$ is the 
comorphism $h : X \dasharrow Z $ over $A\to C$, given by $h_{a}(z):= 
f_{a}(g_{\alpha(a)}(z))$ for $z\in Z_{\beta (\alpha (a))}$. 

Note that for $Y\to B$ and $\alpha :A \to B$, $\alpha^{*}(Y)$ is 
given by the 
$A$-indexed family $\{Y_{\alpha (a)}\mid a \in A\}$.

These set-theoretic descriptions do not depend on cleavages; on the 
contrary, suitably interpreted, reading $a\in A$ etc.\ as  
generalized elements (as in \cite{[SDG]}), they 
describe the universal properties characterizing the objects or maps in 
question  (even in more general categories). Similarly when reading  objects in fibered categories as 
``generalized families'' (as in \cite{JW} Chapter I).

\section{Star bundle functors} As in two the previous sections, we 
consider a category $\B$ with pull-backs, so that we have two 
fibrations over $\B$, $\B^{2}$ and its fibrewise dual $(\B^{2})^{*}$. We also 
have the trivial fibration $\B^{1}$ over $\B$. A {\em star bundle 
functor} (terminology from Kol\'{a}\v{r}, Michor and Slov\'{a}k, \cite{[KMS]}) is now defined to be a functor $S$ over 
$\B$ from $\B^{1}$ to $(\B^{2})^{*}$. By the explicit description in the 
previous section, this amounts to the following data: for each $A\in 
\B$, an arrow $\pi_{A}: S_{0}(A)\to A$, and for each $\alpha : A \to B$, a 
pull-back diagram like (\ref{comorx}), 
$$\begin{diagram}\alpha^{*}S_{0}(B)&\rTo&S_{0}(B)\\
\dTo^{v}& &&\\
S_{0}(A)&  &\\
\dTo^{\pi_{A}} && \dTo_{\pi_{B}} \\
A&\rTo_{\alpha}&B.
\end{diagram}$$

More generally, if $\X\to \B$ is any fibration, a star-bundle functor 
with values in $\X \to \B$ is a functor over $\B$ from $\B^{1}$ to 
$\X^{*}$.

The formation of cotangent bundles  for manifolds is an 
example, to be described in the following Section. It is defined   as the ``fibrewise linear dual'' of the tangent 
bundle, viewed as a vector bundle, i.e.\ as the composite of  $T$ 
with a ``fiberwise duality''  functor $\dagger$, whose categorical 
status will be described.

\section{Vector bundles, and the cotangent bundle}\label{VBCBX}
The full generality of the present Section is probably that of fibrewise 
symmetric monoidal closed category, in the sense of \cite{[Bunge]}, 
\cite{[Shulman]} et al., but we formulate things more concretely  in terms of the fibered category 
$\W \to \B$ of  vector bundles 
(over spaces in a suitable category $\B$ of, say, smooth manifolds). 
Thus $\W _A$ is the category of vector space objects in the category 
$\B /A$. This $\W$ comes with a 
forgetful functor over $\B$ from $\W$ to $\B^{2}$.

Again, we assume a cleavage, and the resulting notation like 
$\alpha^{*}(Y)$ for the chosen pull-back of a vector bundle $Y$ along a 
smooth map $\alpha$. We intend here to clarify the role of the 
notion of fibrewise linear dual of a vector bundle $X\to A$, which we denote 
$X^{\dagger}$ (refraining from using $X^{*}$, since the $*$ already 
has two meanings: $\alpha^{*}$ for pull-backs alomg $\alpha$, and 
$\Y^{*}$ for the fibrewise dual fibration of a fibration $\Y\to \B$). Clearly, this 
dualization is a contravariant endofunctor on $\W_{A}$, for each 
$A\in \B$. (For vector {\em spaces}, the functor $\dagger$ is the 
standard contravariant dualization functor for 
vector spaces.)

\begin{prop}\label{daggerx}The fibrewise linear dualization functor $(-)^{\dagger}: \W_{A} \to 
(\W_{A})^{op}$ extends canonically to a functor 
$\W \to \W^{*}$ over $\B$; it is a Cartesian functor. 
\end{prop}
{\bf Proof.} Consider an arrow over $\alpha$, meaning a commutative diagram
$$\begin{diagram} X&\rTo^{t}&Y\\
\dTo &&\dTo \\
A&\rTo_{\alpha}&B.\end{diagram}$$
with $t$ fibrewise linear, so for each $a\in A$, the map 
$t$ gives a linear map $t_{a}: X_{a}\to Y_{\alpha (a)}$,  
hence a linear  $t_{a}^{\dagger}:(Y_{ \alpha (a)})^{\dagger}\to 
X_{a}^{\dagger}$. But $(Y_{ \alpha (a)})^{\dagger}=(Y^{\dagger})_{ 
\alpha (a)}$. Jointly, these $t_{a}^{\dagger}$ produce a map $ 
\alpha^{*}(Y^{\dagger}) 
\to X^{\dagger}$ in $\W_{A}$, which is the vertical part of the desired 
comorphism; the horizontal part is the arrow $\alpha^{*}(Y^{\dagger}) 
\to Y^{\dagger}$ in 
the diagram defining $\alpha^{*}(Y^{\dagger})$.
 -- If the given square is a pull-back, each $t_{a}$ is 
an isomorphism, hence so is $t_{a}^{\dagger}$, so in this case, the 
vertical part described is an isomorphism; therefore the comorphism 
described is Cartesian; this proves the last assertion.

\medskip

From this perspective, the cotangent bundle construction is a 
functor (over $\B$), namely the composite of the two functors 
$$\begin{diagram}
\B&\rTo^{T}&\W & \rTo^{\dagger}& \W ^{*}
\end{diagram};$$
both $T$ and $\dagger$ are functors over $\B$, hence so is the 
composite. Here, $\W$ and $\W^{*}$ come with forgetful functors to 
$\B^{2}$ and $(\B^{2})^{*}$, respectively. Composing with the 
forgetful functor $\W^{*}\to  (\B^{2})^{*}$  
then gives a functor over $\B$, $\B^{1} \to (\B^{2})^{*}$, i.e.\ a 
star-bundle functor with values in $\B^{2}$.

For the cotangent bundle construction in differential geometry, we
  provide in Section \ref{BV1x} an alternative description of this star-bundle functor, 
without passing through the tangent 
bundle formation, and without using linear structure, namely in terms 
of the ``jet bundle construction'' available in differential geometry.

\medskip

The same argument as for the Proposition gives that
the fibrewise linear dualization functor $(-)^{\dagger}: \W_{A}^{op} \to 
\W_{A}$ extends canonically to a functor 
$\W ^{*}\to \W$ over $\B$; it is likewise Cartesian. (This does not 
depend on whether $V\to V^{\dagger \dagger}$ is an isomorphism.)

\section{Strength}Let $\B$ be a category with finite limits and let $ \X \to 
\B$ and $\Y \to \B$ be  fibrations. 
We consider a functor $F: 
\X \to \Y$ over $\B$ ($F$ is not assumed Cartesian). Then we shall consider a 
certain kind of structure on such a functor, which we call {\em 
fibrational strength}, or just {\em strength}.

For this, we introduce some notation. If $Q\in \B$ and $X\in \X 
_{M}$, we have an object $p^{*}(X) \in \X_{Q\times M}$, where 
$p:Q\times M \to M$ denotes the projection. This object in 
$\X_{Q\times M}$ we denote 
$Q\otimes X$. It comes equipped with a (Cartesian) morphism $Q\otimes 
X \to X$ over $p$.

\medskip

\noindent {\bf Example 1.} Let $\B$ be a category with finite limits. 
Then the codomain fibration $\partial_{1}: \B^{2}\to \B$ is a 
fibration, and  pull-back squares in $\B$ are the Cartesian arrows in 
$\B^{2}$. If $\xi \in \B^{2}$, say $\xi: X \to M$, and $Q \in \B$, 
it is clear that $Q\otimes \xi$, as a map in $\B$, is just $Q\times 
\xi : Q\times X \to Q\times M$.

\medskip  

Let $F: \X \to \Y$ be a functor over $\B$, a {\em 
strength} on $F$ consists in the following data: for $Q\in \B$ and $X 
\in \X_{M}$, one gives a morphism in $\Y _{Q\times M}$
$$t_{Q,M}: Q\otimes F(X) \to F(Q\otimes X),$$
natural in $Q$ and $X$, and satisfying a unit- and associativity 
constraint. 

Since $F$ is a functor over $\B$, there is a canonical  
morphism $F(p^{*}(X))\to p^{*}(F(X))$, i.e.\ $F(Q\otimes X) \to 
Q\otimes F(X)$  
 in $\Y_{M}$. It is invertible if $f$ is Cartesian; so for $F$ Cartesian,
 the  inverse will be    a strength structure  $t_{Q,X}$ on 
$F$. 

If $F: \X \to \Y$ and $G: \Y \to \Z$ are functors over $\B$, equipped 
with strengths $t$ and $s$, respectively, one constructs out of $t$ 
and $s$
 in an evident 
way a strength on the composite functor  $G\circ F: \X \to \Z$.  We obtain
 a 2-category: objects are categories 
fibered over $\B$, arrows are functors over $\B$ equipped with 
strength, and 2-cells are the vertical natural transformations  between parallel 
functors over $\B$,  compatible with the given strengths.

If $\B$ is the category of smooth manifolds, a map $h:Q\times M \to 
N$ in $\B$ may be seen as a smoothly parametrized family of smooth 
maps $h(q,-): M 
\to N$ (with $Q$ as the space of parameters), and a map $H: Q\otimes X 
\to X'$ over $h$ may be seen as a $Q$-parametrized  family $H$ of maps in $\X$ 
from $X$ to $X'$, with the $q$th member of this family living over 
$h(q,-):M \to N$. 

A strength $t$ of $F: \X \to \Y$ gives rise to a process transforming a 
paramterized family of maps in $\X$ to a similarly parametrized 
family of maps in $\Y$, as follows. Given a map $h:Q\times M \to N$ 
in $\B$. Then given a map in $\X$ over $h$, say  $H: Q\otimes X \to 
X'$, then 
the composite
$$\begin{diagram}Q\otimes F(X)& \rTo^{t_{Q,X}}& F(Q\otimes X) &\rTo 
^{F(H)}&F(X')
\end{diagram}$$
is  a map in $\Y$ over $h$; so  $F$ has has transformed the 
$Q$-parametrized family of maps from $X$ to $X'$ into a 
$Q$-parametrized family of maps $F(X) \to F(X')$.  This property of 
$F$ is called {\em regularity} in \cite{[KMS]}.

\medskip

The identity functor $\B \to \B$ is a fibration,  denoted 
$\B^{1}$; for this fibration 
$Q\otimes X = Q\times X$. Recall that a {\em bundle  functor} is a 
functor $F: \B^{1} 
\to \B^{2}$ over $\B$, thus for $X \in \B$, the object $F(X)$ in 
$\B^{2}$ is an arrow of the form $F(X): F_{0}(X) \to X$ in $\B$, 
where $F_{0}$ is the composite of $F$ with the domain formation 
$\partial _{0}:\B^{2}\to \B$. A 
fibrational strength $t$ of such functor gives for $Q$ and $X$ in $M$ 
an arrow $t_{Q,X}$ in $\B^{2}$ from $Q\otimes F(X)$ to $F(Q\times X)$, which 
amounts to a commutative  square in $\B$ (really just a triangle) of the form
\begin{equation}\label{tdoubx}\begin{diagram}Q\times F_{0}(X)&\rTo^{t''_{Q,M}} & F_{0}(Q\times X)\\
\dTo^{Q\otimes F(X)}&t_{Q,X}&\dTo_{F(Q\times X)}\\
Q\times X&\rTo _{\id}&Q\times X
\end{diagram},\end{equation}
and the top map in this square (as $Q$ and $X$ range over $\B$) 
equips the endofunctor $F_{0}: \B \to \B$ with a {\em tensorial} 
strength $t''$ in the sense of \cite{[SFMM]}. Vice versa,  if such $t''$ 
make the squares like the above commute, these squares will constitute 
a fibrational strength $t$ on $F$. (To say that the squares commute is in 
turn equivalent to saying that $F$, viewed as a natural transformation 
from $F_{0}$ to the identity functor on $\B$, is a strong natural 
transformation, in the sense of tensorial strength.)

\medskip

Let $F:\B^{1}\to \B^{2}$ be a bundle functor preserving finite 
products. Thus  $F(Q\times B)\cong F(Q)\times F(B)$ by the canonical 
map. In particular (since 
$\partial _{0}: \B^{2}\to \B$ preserves products), $F_{0}(Q\times B) 
\cong F_{0}(Q)\times F_{0}(B)$. An example is where $F$ is the tangent 
bundle functor (ignoring the fibrewise linear structure).

\medskip
A particular bundle functor on $\B$ is the diagonal $\Delta$ 
associating to $B \in B$ the identity map $B \to B$. It terminal 
among bundle functors $\B^{1}\to \B^{2}$.

A {\em section} of a bundle functor $F$ is a natural transformation 
$z$ 
(over $\B$) 
from  $\Delta :\B^{1}\to \B^{2}$ to $F$, thus to 
each $B\in \B$, $z_{B}:B \to F_{0}(B)$ is a section of 
$F(B):F_{0}(B) \to B$. The zero section of a tangent bundle is an 
example.

\begin{prop} Let $F$ be a finite-product preserving bundle functor 
equipped with a zero section. Then $F$ 
carries a canonical (fibrational) strength.
\end{prop}
{\bf Proof/construction.} By the above (cf.\ (\ref{tdoubx}), it 
suffices to   construct in $\B$ a map 
$t''_{Q,B}: Q\times F_{0}(B) \to F_{0}(Q\times B)$. This is taken to be 
the composite of $z_{Q}\times F_{0}(B): Q\times F_{0}(B) \to 
F_{0}(Q)\times F_{0}(B)$ with the isomorphism $F_{0}(Q)\times 
F_{0}(B)\cong F_{0}(Q\times B)$.

\medskip

The notion of (fibrational) strength of a functor $F: \X \to \Y$ over $\B$, 
in the sense described here, generalizes the 
notion of ``regularity'' of a bundle functor, \cite{[KMS]} 14.21 (and 
18.10). 
The reason we change terminology from ``regularity'' to ``strength'' 
is to emphasize  1) that, in the abstract setting, it is a {\em structure}
 on the functor in 
question, not just a property, and 2) to tie it up with the notion of 
(monoidal, or tensorial) strength considered in the context of endofunctors $F$ on a 
monoidal category $\B$, as in    
\cite{[SFMM]}, 
(or  \cite{[Moggi]}, 
or \cite{[CMTD]} Section 2, for a recent account).
Such a structure  in turn is equivalent to a 
$\B$-enrichment of $F$, in case $\B$ is monoidal closed; cf.\ \cite{[SFMM]}.

\medskip 

If $\X\to \B$ is the fibration, where $\X_{B}$ is the category of 
vector space objects in $\B/B$ (or group objects, or any other 
algebraic kind of structure), then there is a faithful forgetful functor $\X 
\to \B^{2}$ over $\B$, which is Cartesian, in particular, it 
preserves the formation $Q\otimes X$. So if the bundle functor $F$ 
considered above factors through $\X \to \B^{2}$ (a ``vector bundle 
functor''), then the $t''_{Q,B}$ constructed above is the underlying 
arrow of an arrow in $\X$, i.e.\ is fibrewise linear, and equippes 
the vector bundle functor $F$ with a fibrational strength.

This in particular applies to the tangent bundle formation. The 
cotangent bundle functor likewise carries a canonical strength, by 
the following
\begin{prop}Let $\X \to \B$ be the category of vector bundles. Then 
if there is given a strength on $F: \B \to \X$, then there is 
canonically associated a strength on the star bundle functor 
$F^{\dagger}: 
\B \to \X^{*}$.\end{prop}
{\bf Proof.} This follows since the dualization functor $\dagger : \X 
\to \X^{*}$ is Cartesian, 
and hence carries a canonical strength; and a composite of two 
functors with a strength has a strength. Note that the instantiations 
$t_{Q,B}$ of the strength described here are (vertical) maps $Q\otimes 
F^{\dagger}(B) \to 
F^{\dagger}(Q\times B)$ in $\X^{*}$, and hence as vector bundle maps (maps in $\X$) are maps
$F^{\dagger}(Q\times B) \to Q\otimes F^{\dagger}B$; for,  the fibre 
of $\X^{*}_{Q\times B}$ is $(\X_{Q\times B})^{op}$.

\section{Flow natural transformation}

Consider a category $\B$ with Cartesian products, and consider an 
endofunctor $F: \B \to \B$ with a tensorial strenght $t''_{Q,B}:Q\times 
F(B) \to F(Q\times B)$ (for all $Q$ and $B$ in $\B$). If now $D\in \B$ is an 
exponentiable object, so $(D\times - )\dashv (D\pit -)$, one derives, 
for all $B\in \B ,$ 
a map $$\lambda_{D,B}:F(D\pit B) \to D\pit F(B),$$ namely the 
exponential transpose of the composite
$$\begin{diagram}D \times F(D\pit B)&\rTo^{t''_{D,D\pit B}}&F(D \times 
(D\pit B))&\rTo^{F(\ev)}&F(B).$$ 
\end{diagram}.$$
If $D_{1} \to D_{2}$ is a map between exponentiable objects, one gets a map 
$D_{2}\pit B \to D_{1}\pit B$ (`` $\pit$ is contravariant in the first 
variable "), and then $\lambda$ will be natural in the $D_{i}$'s, in 
an evident sense. Also, $\lambda _{1,B}: F(1\pit B)\to 1\pit F(B)$ 
may be identified with the identity map on $F(B)$. So if the 
exponentiable object $D$ is 
equipped with a point $0:1\to D$, one obtains a commutative triangle
\begin{equation}\label{lamx}\pi\circ \lambda_{D,B}= F(\pi 
),\end{equation} where 
$\pi$ denotes $0\pit B$ or $0\pit F(B)$. So $\lambda_{D,B}$ is a map 
of bundles over $F(B)$.

If $\B $ is Cartesian closed (or even just symmetric monoidal closed), all objects $D$ are exponentiable, 
and therefore one 
has such a $\lambda_{D,X}$ for all $D,X$, and this data encodes the 
strength of $F$ in what may be called the  ``cotensorial'' form of
strength, cf.\ \cite{[CCGBCM]} or \cite{[CMTD]}.

\medskip

For the following, we shall assume that $\B$ is a model for synthetic 
differential geometry, in 
particular, 
it contains the category of smooth manifolds, but also it contains 
some ``infinitesimal objects'', in particular, it contains an object 
$D$  with the property that for any manifold $M$, $T(M)= D\pit M$, and 
the base map $T(M) \to M$ is ``evaluation at $0\in D$'', where 
$0:{\bf 1}\to D$ is a given point of $D$. Then for any  endofunctor $F: \B \to 
\B$ with a  (tensorial) strength,  we have, by the above 
construction, 
$\lambda _{D,M}: F(D\pit M) \to D\pit F(M)$. 
If $F(M)$ is a manifold whenever $M$ is, then this map is a map 
between manifolds, since manifolds form a full subcategory of $\B$;
and $\lambda_{D,M}$ is natural in $M$ (since $\lambda$ is); 
it is the {\em flow natural transformation} $F(T(M))\to T(F(M))$ for $F$ considered in 
\cite{[KMS]} 39.1 (denoted there $\iota _{M}$). 
It originated in a discussion between Kol\'{a}\v{r} and the present 
author in the early 1980s, see the ``Remarks'' at p.\ 349 in loc.cit.

An application of the flow natural transformation is that it gives a 
``prolongation procedure'' for vector fields on $M$: to a vector 
field $\xi : M \to T(M)$ on $M$, one constructs a vector field 
$\tilde{\xi}:F(M) \to T(F(M))$ 
on $F(M)$, namely the composite
$$\begin{diagram}F(M)&\rTo^{F(\xi )}&F(T(M))& \rTo^{\lambda 
_{D,M}}&T(F(M)).
\end{diagram}$$

\section{Jet bundles}\label{JBX} The $k$th order jet bundle  of a smooth fibre 
bundle $p:E\to B$ in differential geometry is another smooth fibre 
bundle $J^{k}(p)\to B$ (usually just denoted  $J^{k}(E)$). The fibre 
over $b\in B$ consist of $k$-jets at $b$ of sections of $E\to B$, see 
\cite{[PALAIS]}, or \cite{[SGM]} 2.7 (and Remark 7.3.1); in the latter synthetic context, 
the notion of $k$-jet becomes representable, in the sense that 
there is for every $b\in B$ a subset $\FM _{k}(b) \subseteq B$ 
(with $b\in \FM_{k}(b)$), such 
that a $k$-jet at $b$ is a map with domain $\FM_{k}(b)$, in 
particular, a $k$-jet of a section of $p:E\to B$ is a  
map $s: \FM_{k}(b) \to E$ with $p\circ s$ equal to the inclusion map 
$\FM_{k}(b) \to B$.  For fixed $B$, $J^{k}(E)$ depends in a functorial 
way on $E$ in the category of smooth fibre bundles over $B$; 
synthetically, if $f:E \to E'$ is a map of bundles over $B$ and $s$ 
is a $k$-jet of a section of $E$, $J^{k}(f)$ takes a section 
$s:\FM_{k}(b) \to E$ of 
$E$ to the map $f\circ s: \FM_{k}(b) \to E'$. Or, in classical set-up, 
post-composition by $f$ of a partial section $s$, representing the 
given jet, has as $k$ jet at $b$  the desired jet section of $E'$.

So for each $B\in \B $, we have an endofunctor $J^{k}$ on the category $\B  
_{B}$ (= the category of smooth fibre bundles over $B$). We shall 
investigate the functorality properties of $J^{k}$ as $B$ varies over 
the category of smooth manifolds. 

To simplify the exposition, we shall embed the category  of 
smooth manifolds into a topos model $\B $ of synthetic differential 
geometry (cf.\ e.g.\ \cite{[SDG]} or \cite{[SGM]}), where the jet 
construction works not just for smooth fibre bundles $E\to B$ 
but for any smooth map $E\to B$, (where $B$ is a manifold), so that for each $B$, $J^{k}$ is an 
endofunctor $J^{k}_{B}$ on $\B  /B$. 

\medskip
 The description (from \cite{[SGM]} Remark 7.3.1)
 of $J^{k}$ is given in terms of the locally Cartesian closed 
structure of $\B$, as follows: The data of the $\FM_{k}(b)$, as $b$ ranges over 
$B$, resides in ``the $k$th neighbourhood of the 
diagonal\footnote{The use of a ``$k$th  neighbourhood of 
the diagonal'', also called ``prolongation spaces'', for the consideration of jet bundles is crucial in 
\cite{KumpS}; the setting there is that of  manifolds equipped with a 
structure sheaf of rings (that may contain nilpotent elements), as 
considered by Grothendieck and Malgrange. }of $B$'',
$$\begin{diagram}B_{(k)}&\pile{\rTo^{c}\\ \rTo_{d}}& B;
\end{diagram}$$
and similarly for $A$. Here, $B_{(k)}\subseteq B\times B$ consists of 
pairs $(b,b')$ with $b'\in \FM_{k}(b)$, and $c$ and $d$ are the 
restrictions of the 
two projections $B\times B \to B$. Similarly for $A_{(k)}\subseteq 
A\times A$ (where we  again denote the two projections 
by $c$ and $d$). The map $\alpha \times \alpha : A\times A \to B \times B$ 
restricts to a map $\overline{\alpha}: A_{(k)}\to B_{(k)}$ 
(equivalently, any map $A\to B$ restricts to a map $\FM_{k}(a)\to 
\FM_{k}(\alpha (a))$). Pulling 
back along $d:B_{(k)}\to B$ defines a functor $d^{*}:\B/B \to 
\B/B_{(k)}$, and since $\B$ is locally Cartesian closed, this functor 
has a right adjoint $\Pi_{d}: \B/B_{(k)} \to \B/B$. In these terms, 
the endofunctor $J^{k}$ on $\B/B$ is just the composite $\Pi_{d}\circ 
c^{*}$.

\begin{thm}The functors $(J^{k}_{B} )^{op}$ are the fibres of an 
endofunctor $J^{k} : (\B ^{2})^{*} \to (\B ^{2})^{*}$ over 
$\B $.
\end{thm}
{\bf Proof.} Since we already have the functor $J^{k}$ on the individual 
$\B/A$, (for $A\in \B$) it is possible to prove this using the construction 
Theorem \ref{constrx}; however, since the categories and functors have so concrete 
descriptions, it is easier and more informative to give the 
construction and proofs ad hoc, using set/family - theoretic 
descriptions, as in the Remark at the end of Section \ref{FDX}.
The construction amounts to a process which to a comorphism $f$ over 
$\alpha:A\to B$ from $X\to A $ to $Y\to B$ associates a comorphism 
$J^{k}f$
over $\alpha$ from $J^{k}X$ to $J^{k}Y$. 
Recall that in the set theoretic description (translating 
(\ref{comorx}) into elementwise terms), a comorphism $f$ over 
$\alpha$ amounts to a family 
of maps $f_{a}:Y_{\alpha(a)}\to X_{a}$, for $a$ ranging over $A$. 
Similarly, the required $J^{k}f$ is to consist of a family 
$(J^{k}f)_{a}: (J^{k}Y)_{\alpha(a)}\to (J^{k}X)_{a}$. An element 
$s$ in $(J^{k}Y)_{\alpha(a)}$ is a partial section $s: \FM_{k}(\alpha (a)) 
\to Y$ of $Y\to B$. The composite
$$\begin{diagram}\FM_{k}(a)&\rTo^{\alpha}& \FM_{k}¥(\alpha (a))&\rTo^{s}&Y
\end{diagram}$$
is a map $\FM_{k}(a) \to Y$ over $\alpha$, or, equivalently, a map 
$\FM_{k}(a) \to \alpha^{*}(Y)$ over $A$, 
thus an element of $(J^{k}(\alpha^{*}Y))_{a}$, to which we may apply 
the map  
$J^{k}f : J^{k}\alpha^{*}(Y) \to J^{k}X$; this means just: 
post-composing with $f: \alpha^{*}Y \to X$. Thus, the element in 
$(J^{k}X)_{a}$ that we get is the map $\FM_{k}(a) \to X$ given 
elementwise as follows: to input $a'\in \FM_{k}(a)$,  we get as output
$f(s(\alpha (a')))\in X_{a'}$. From this later 
description, the compatibility of the construction of $J^{k}$ with 
composition of 
comorphisms is almost immediate.

\medskip

\noindent{\bf Remark 1.} Let us note that if the fibres of $Y \to B$ carry some algebraic 
structure, say that of vector spaces, then so do the fibres of 
$J^{k}Y$. This follows, since $J^{k}=\Pi_{d}\circ c^{*}$ is a 
composite of two right adjoints, so preserves algebraic structure. So $J^{k}: (\B^{2})^{*} \to (\B^{2})^{*}$ 
lifts to a 
functor $\W^{*}\to \W^{*}$, where $\W \to \B$ is the category of 
vector bundles. Similarly for other kinds of algebraic structure, 
e.g.\ {\em pointed} spaces.

\medskip

\noindent{\bf Remark 2.} Let us also remark that the existence of the maps $\alpha^{*}(J^{k}Y)
\to J^{k}(\alpha^{*}Y)$ considered above implies the existence of a 
fibrational strength of the functor $J^{k}$: just take  $\alpha $ 
to be the projection $Q\times B \to B$.

\section{Bundle valued 1-forms}\label{BV1x}
The natural setting for the present subsection is the fibration of 
{\em pointed} bundles (with morphisms preserving the given points); 
there is a  forgetful functor from the fibration $\W \to \B$ of vector 
bundles to the codomain fibration  $\B^{2}\to \B$, and this functor factors through 
the fibration of pointed bundles, but in order not to overload the exposition with too 
much terminology and notions, the presentation is in terms of the 
fibration $\W \to \B$ of vector bundles. If $E\to B$ is such a 
bundle, a 1-jet at $b\in B$ of a section, i.e.\ a partial section 
  $s: \FM _{1}(b) \to E$ is called an {\em $E$-valued} 
(cominatorial) {\em  cotangent} at 
$b$ if $s(b)=0_{b}$. So $s$ and the zero section agree on $b\in 
\FM_{1}(b)$, but do not necessarily agree on the whole of $\FM 
_{1}(b)$. Clearly, the set of $E$-valued cotangents form a sub-bundle 
of $J^{1}(E)\to B$, called the {\em bundle of $E$-valued 1-forms}; 
let us denote it $\Omega ^{1} (E) \to B$. It is a subfunctor  of the functor
$J^{1}: \W^{*} \to \W^{*}$.

Let $R$ be a fixed vector space (typically, the ground field). There is a functor 
$\B \to \W$ over $\B$, assigning to $B\in \B$ the constant vector bundle 
$B\times R \to B$. This functor is Cartesian, and hence may equally well be viewed as 
a functor $\B\to \W^{*}$, since the category of Cartesian 
arrows in $\W$ and $\W^{*}$ are the same.
Composing with $\Omega ^{1}$, 
$$\begin{diagram}\B &\rTo ^{-\times R}&\W^{*}& \rTo ^{\Omega ^{1} }& \W^{*}
\end{diagram},$$
we get a vector star bundle functor, i.e.\ a functor  $\B \to 
\W^{*}$; this is (isomorphic to) the cotangent bundle functor 
$T^{\dagger}: \B \to \W^{*}$ described Section \ref{VBCBX}.
In fact, the  bundle functor $T:\B \to \W$ (tangent bundle)  is sometimes {\em defined}
as the composite 
$$\begin{diagram}\B&\rTo^{\Omega } &\W^{*}& 
\rTo^{\dagger}&\W.
\end{diagram}$$

\small
\noindent Anders Kock\\
Dept.\ of Math., University of Aarhus, Denmark\\
e-mail: kock at math dot au dot dk\\
May 2015


\begin{thebibliography}{99}
\bibitem{[BGN]} C.\  Barwick, S.\ Glasman and D.\ Nardin, Dualizing Cartesian and 
Cocartesian Fibrations, arXiv:1409.2165v1 [math.CT], 2014.

\bibitem{[Borceux]} F.\ Borceux, {\em Handbook of Categorical Algebra}  Vol.\ 2, Encyclopedia of 
Math.\ and its Applications 51,   Cambridge 1994.

\bibitem{[Bunge]} M.\ Bunge, Tightly bounded completions, Theory and 
Applications of Categories Vol.\ 28, No.\ 8, 
213-240, 2013.

\bibitem{[Elephant]} P.\ Johnstone, {\em Sketches of an Elephant; A Topos Theory 
Compendium}, Oxford Logic Guides 43-44, Oxford 2002. 

\bibitem{JW} P.\ Johnstone and G.\ Wraith, Algebraic theories in 
toposes, in {\em Indexed Categories and Their Applications} (ed.\ 
P.\ Johnstone and R.\ Par\'{e}), Springer Lecture Notes in Math.\ 661 
1978.







\bibitem{[CCGBCM]} A.\ Kock, Closed categories generated by commutatice monads, 
J.\ Austral.\ Math.\ Soc.\ 12 , 405-424, 1971. 

\bibitem{[SFMM]} A.\ Kock, Strong functors and monoidal monads, Arch.\ Math.\ 
(Basel) 23, 113-120 (1972).



\bibitem{[SDG]}  A.\ Kock, {\em Synthetic Differential Geometry}, London Math.\ Soc.\ 
Lecture Notes 51, Cambridge  1981 (2nd ed.\ London Math.\ Soc.\ 
Lecture Notes 333, Cambridge 2006).

\bibitem{[SGM]} A.\ Kock, {\em Synthetic Geometry of Manifolds}, Cambridge Tracts in 
Mathematics 180, Cambridge 2010.


\bibitem{[CMTD]} A.\ Kock, Commutative monads as a theory of distributions,  Theory and 
Applications of Categories Vol.\ 26, No.\ 4, 2012. 


\bibitem{[DFET]} A.\ Kock, The dual fibration in elementary terms, 
arXiv:1501.01947v1 [math.CT], 2015.

\bibitem{[KMS]} I.\ Kol\'{a}\v{r}, P.\ Michor and J.\ Slov\'{a}k, 
{\em Natural 
Operations in Differential Geometry}, Springer 1993.

\bibitem{KumpS}A.\ Kumpera and D.\ Spencer, {\em Lie Equations}, Annals of 
Mathematics Studies 73, 1972.

\bibitem{[Moggi]} E.\ Moggi, Notion of computations and monads, Information and 
Computation 93, 55-92, 1991.


\bibitem{[PALAIS]} R.\ Palais, {\em Seminar on the Atiyah-Singer Index Theorem}, 
Annals of Mathematics Studies 57,   1965.

\bibitem{[Shulman]} M.\ Shulman, Enriched indexed categories, Theory and 
Applications of Categories Vol.\ 28, No.\ 21, 
616-695, 2013.

\bibitem{[Weber]} M.\ Weber, Polynomials in categories with pullbacks, 
Theory and 
Applications of Categories Vol.\ 30, No.\ 16, 
533-598, 2015.

\end{thebibliography}
 \end{document}